\documentclass[12pt]{article}
\usepackage{amsmath}
\usepackage[russian,english]{babel}
\usepackage{amsthm}
\usepackage{epsfig}
\input epsf
\def\Ima{\operatorname{Im}}
\def\Rea{\operatorname{Re}}
\def\sgn{\operatorname{sgn}}
\def\card{\operatorname{card}}
\def\sect{\operatorname{sect}}
\def\tr{\operatorname{tr}}
\def\Ima{\operatorname{Im}}
\def\R{\mathbf{R}}
\def\C{\mathbf{C}}
\def\H{\mathbf{H}}
\def\bC{\mathbf{\overline{C}}}
\def\const{\mathrm{const}}
\def\Conf{\mathrm{Conf}}
\def\Sub{\mathrm{Sub}}
\def\Combs{\mathrm{Combs}}
\newtheorem{theorem}{Theorem}

\newtheorem{lemma}{Lemma}
\theoremstyle{remark}
\begin{document}
\author{
Alexandre Eremenko\thanks{Supported by NSF grant
DMS-1067886.}$\;$ and Peter Yuditskii\thanks{Supported by the Austrian Science Fund FWF, project no: P22025-N18.}}
\title{Comb functions}
\maketitle
\begin{abstract} We discuss a class of regions and conformal mappings
which are useful in several problems of approximation theory,
harmonic analysis and spectral theory.\footnote{
This text was prepared for a plenary talk given by P. Yuditskii
on 11-th International Symposium on Orthogonal Polynomials,
Special Functions and Applications, dedicated to celebrate Francisco
(Paco) Marcell\'an's 60-th birthday.}

MSC: 30C20, 41A10, 47B36, 41A50. Keywords: conformal map, Green function,
Martin function, uniform approximation, Jacobi matrices, Riesz bases,
spectral theory.

\end{abstract}

\section{Introduction}

We begin with two simple classical problems which serve as motivation.
Then in sections 2--5 we describe some classes of regions, corresponding 
conformal maps and entire and subharmonic functions.
In sections 6--7 we discuss various problems where these classes appear. 
\vspace{.1in}

1. Polynomials of least deviation from zero.
Let $E\subset\R$ be a compact set on the real line, and $P_n$ 
a polynomial 
with minimal sup-norm $L_n=\| P_n\|_E$ among all monic 
polynomials of degree $n$.

If $n<\card E$, then $P_n$ is unique and can be characterized by
the following properties:

(i) $P_n$ is real, and all its zeros are real and simple,

(ii) For every pair of adjacent zeros $x_1<x_2$ there is a point $y\in(x_1,x_2)\cap E$
such that $|P_n(y)|=L_n$.

(iii) At the points $a_1=\inf E$ and $a_2=\sup E$, we have $|P_n(a_j)|=L_n$.

For a simple variational argument which proves (i)--(iii) see \cite{Akh2,SY1}.
These polynomials $P_n$ can be represented in terms of 
special conformal maps.

Let $m,k$ be integers, $k-m=n$, and let $D$ be a region obtained from the half-strip
$$\{ z=x+iy: \pi m<x<\pi k,\; y>0\}$$
by removing vertical intervals $\{\pi j+it:0\leq t\leq h_j^\prime\}$, $m<j< k$,
where $h_j^\prime\geq 0$, see Fig. 1 (right).
\begin{center}
\epsfxsize=5.0in
\centerline{\epsffile{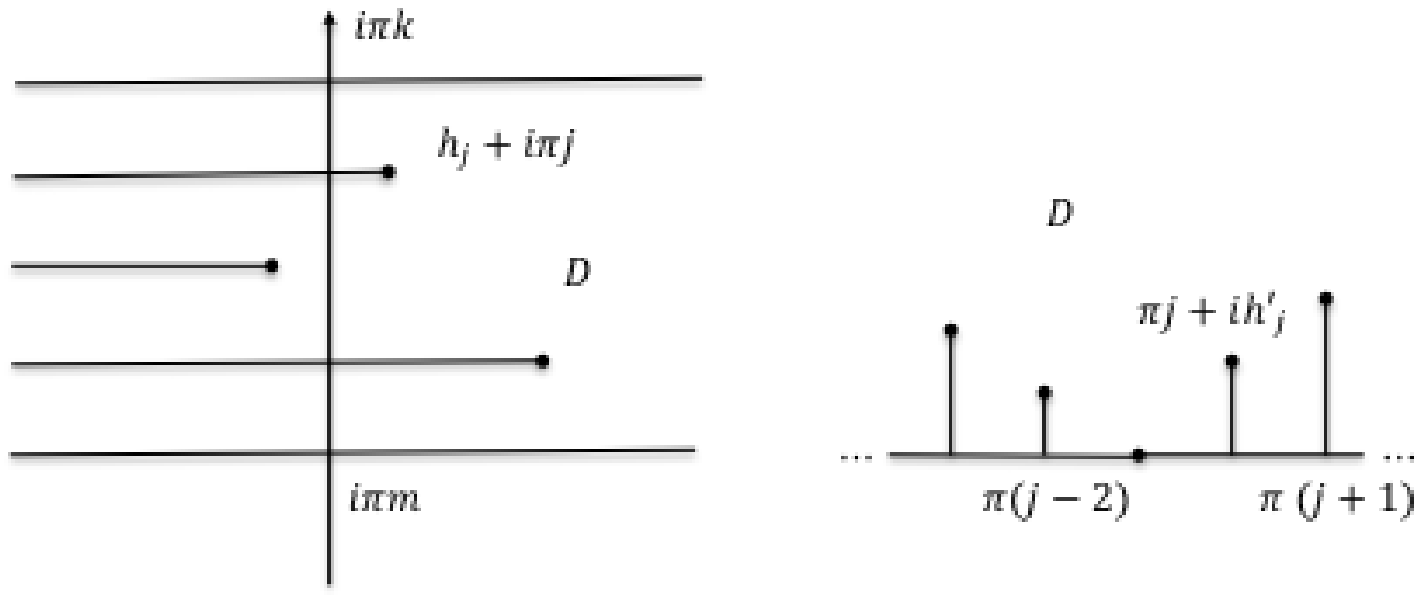}}
Fig. 1. Comb regions of V (left) and MO (right) types.
\end{center}
Let $\theta$ be a conformal map from the upper half-plane $\H$ to $D$,
such that $\theta(\infty)=\infty,\; \theta(a_1)=\pi m,\;\theta(a_2)=\pi k.$
Then $z\mapsto\cos\theta(z)$ is an analytic function in the upper half-plane,
which is real on the real line. So it extends to an entire function, and
the behavior at $\infty$ shows that this entire function is a polynomial
of degree $n$. Choose $L$ so that the polynomial $P=L\cos\theta$ is monic.
It is easy to check that our extremal polynomials satisfying
(i)--(iii) are of this form,
with an appropriate choice of parameters $h_j^\prime\geq 0$.

The set $E$ is contained in $E':=\theta^{-1}([\pi m,\pi k])$. 
This set $E'$ is the maximal extension of $E$, for which the extremal
polynomial is the same as the one for $E$.
Critical points of $P_n$ are preimages of the tips of the slits under $\theta$,
and critical values are $\pm\cosh h_j^\prime.$ The $\theta$-preimages
of the points $\pi j,\; m<j<k$ are solutions of $P_n(z)=\pm L$, and all
these solutions are real.

For example, if $E=[-1,1]$ we take all $h_j^\prime=0$, and $P_n$ is the $n$-th
Chebyshev polynomial. If $E$ consists of two intervals symmetric
with respect to $0$, and $n$ is even, we take all $h_k=0$, except one,
$h_{(m+k)/2}^\prime>0$. 

On polynomials of least deviation from $0$ on several intervals we refer to
\cite{Akh1,Akh2}, \cite[vol. 1]{Akh3} and the survey \cite{SY1},
where the representation
$P_n=L\cos\theta$ is used systematically.
\vspace{.1in}

2. Spectra of periodic Jacobi matrices.
Consider a doubly infinite, periodic Jacobi matrix
$$J=\left(\begin{array}{cccccc}
\ldots&\ldots&\ldots&\ldots&\ldots&\ldots\\
                               p_{-1}&q_{-1}&p_{0}&0&0&0\\
                                0    &p_{0}&q_0&p_1&0&0\\
                                0    &  0   &p_1&q_1&p_2&0\\
                                0    &  0   & 0 & p_2&q_2&p_3\\
\ldots&\ldots&\ldots&\ldots&\ldots&\ldots\end{array}\right)$$
which is constructed of two periodic sequences of period $n$, where $q_j$
are real, and $p_j>0$. This matrix defines a bounded self-adjoint operator
on $\ell_2$, and we wish to describe its spectrum \cite{GK,March1,Simon}.

To do this we consider a generalized eigenvector $u\in\ell_\infty$ which
satisfies
$$Ju=zu,\quad z\in\C.$$
For fixed $z$, this can be rewritten as a recurrent relation on the
coordinates of $u$:
$$p_{j+1}u_j+(q_{j+1}-z)u_{j+1}+p_{j+2}u_{j+2}=0,$$
which we rewrite in the matrix form as
$$\left(\begin{array}{c}u_{j+1}\\ p_{j+2}u_{j+2}\end{array}\right)=
\left(\begin{array}{cc}0&1/p_{j+1}\\-p_{j+1}&(z-q_{j+1})/p_{j+1}\end{array}\right)
\left(\begin{array}{c}u_j\\ p_{j+1}u_{j+1}\end{array}\right).$$
Thus 
$$\left(\begin{array}{c}u_n\\ p_{n+1}u_{n+1}\end{array}\right)=T_n(z)
\left(\begin{array}{c}u_0\\ p_1u_1\end{array}\right),$$
where $T_n(z)$ is a polynomial matrix with determinant $1$, which is called
the transfer-matrix. To have a bounded generalized eigenvector $u$,
both eigenvalues of $T_n$ must have absolute value $1$.
This happens if and only if
$$|P_n(z)|:=|\tr T_n(z)|/2\leq 1.$$
As $P_n$ is a real polynomial, the spectrum is the preimage of the interval
$[-1,1]$. As our matrix $J$ is symmetric, the spectrum must be real,
this is the same as the condition that all solution of the equations
$P_n(z)=\pm1$ are real,
so we obtain a polynomial of the same kind as in Example 1. 
For every real polynomial with this property,
there exists a periodic Jacobi matrix whose spectrum is $P_n^{-1}([-1,1])$,
and all matrices $J$ with a given spectrum can be explicitly described
\cite{March2,Simon}.
Our polynomial has a representation $P_n=\cos\theta$,
where $\theta$ is a conformal map
of the upper half-plane onto a comb region $D$ as in Example 1.
We obtain the result that the spectrum of a periodic Jacobi matrix consists
of the intervals -- preimage of the real line under a conformal map $\theta$.

We can prescribe an arbitrary sequence $h^\prime_j,\;1\leq j\leq n-1$,
construct a conformal map $\theta:\H\to D$, where $D$ is the region
shown in Fig. 1 (right), and the polynomial $P=\cos\theta$ will have
critical values $(-1)^j\cosh h_j^\prime$ and all solutions of $P(z)=\pm1$
will be real. Such polynomial $P$ is defined by its critical
values of alternating sign up to a change of the independent variable
$z\mapsto az+b,\; a>0, b\in\R$.
Later we will show that any real polynomial with arbitrary
real critical points is defined by its critical
values up to a change of the independent variable
$z\mapsto az+b,\; a>0,b\in\R$.

\section{Comb representation of $LP$ entire functions}

In both examples in the Introduction, the class of real polynomials $P$
such that all solutions
of $P(z)=\pm1$ are real appears. Evidently, all zeros of such polynomials
must be real and simple.
Here we discuss a representation of polynomials with real zeros,
not necessarily simple,
using conformal mappings and generalization of this representation
to a class of entire functions.

Let $P$ be a non-constant real polynomial of degree $n$ with all zeros real.
Let $\varphi=\log P$ be a branch
of the logarithm in the upper half-plane $\H$. Then
$$-\varphi'=-\frac{P'}{P}=-\sum_{n=j}^n\frac{1}{z-z_j}$$
is an analytic function in $\H$ with positive imaginary part.  
\begin{lemma}\label{lemma1} An analytic function $\psi$ in $\H$ whose derivative has positive
imaginary part is univalent.
\end{lemma}
{\em Proof.} Suppose that $\psi(z_1)=\psi(z_2)$, $z_j\in \H,\; z_1\neq z_2.$
Then 
$$0=\frac{\psi(z_1)-\psi(z_2)}{z_1-z_2}=
\int_0^1\psi'(z_2+t(z_1-z_2))dt,$$
but the last integral has positive imaginary part and thus cannot be $0$.
\vspace{.1in}

It is easy to describe the image $\varphi(\H)$. By Rolle's theorem, all
zeros of $P'$ are real and we arrange them in a sequence
$x_1\leq\ldots\leq x_{n-1}$ where each zero is repeated according
to its multiplicity. Let $c_j=P(x_j)$ be the {\em critical sequence}
of $P$. Then the region $D=\varphi(\H)$ is obtained
from a strip by removing $n-1$ rays:
\begin{equation}\label{V-comb}
D=\{ x+iy:\pi m<y<\pi k\}\backslash\bigcup_{m<j<k}
\{ x+i\pi j:-\infty<x\leq h_j\}.
\end{equation}
Here $m-k=n$, and $h_j=\log|c_{k-j}|\geq-\infty.$
Thus 
\begin{equation}\label{vc}
P=\exp\varphi,
\end{equation}
where $\varphi$ is a conformal map of the upper half-plane onto a region $D$
we just described. Such regions will be called polynomial {\em $V$-combs}.
Letter V in this notation is used because this type of representation
was introduced by  Vinberg in \cite{vinb}.

Now suppose that an arbitrary finite sequence $h_j\in [-\infty,\infty),\;
m<j<k$
is given. Consider a $V$-comb $D$ corresponding to this sequence, and
a conformal map $\varphi:\H\to D$.  Using the same argument with reflection
as in Example 1 in the Introduction, it is easy
to see that $\exp\varphi$ is a real polynomial of degree $n$ with critical
values $(-1)^je^{h_j}$. 

We obtain
\begin{theorem} For every finite sequence $c_1,\ldots,c_{n-1}$
with the property 
\begin{equation}\label{alt}
c_{j+1}c_j\leq 0,
\end{equation}
there exists a real polynomial
with real zeros for which this sequence is the critical sequence.
Such polynomial is defined by its critical sequence up to a real
affine change $z\mapsto az+b,\; a>0$ of the independent variable.
\end{theorem}

Now we extend this result to entire functions.
Recall that an entire function belongs to the class $LP$ (Laguerre-P\'olya)
if it is a limit of real polynomials with all zeros real.
For more information on the $LP$-class and its applications we refer
to \cite{L1}.
Consider the following class of regions. 
Begin with 
$\{ x+iy:\pi m<y<\pi k\}$, where $m,k$ are integers or $\pm\infty$,
$-\infty\leq m<k\leq\infty$
and
remove from this region the rays of the form
\begin{equation}\label{rays}
\{ x+i\pi j:x\leq h_j\}, \; m<j<k,\end{equation}
where $h_j\in[-\infty,\infty)$. 
A region of this form is called a $V$-comb corresponding
to a sequence $(h_j)$, $h_j\in[-\infty,\infty)$.
The sequence can be finite, or infinite in one direction, or
infinite in both directions.

\begin{theorem}
The following statements are equivalent:
\vspace{.1in}

\noindent
(i) $f\in LP$,
\vspace{.1in}

\noindent
(ii) $f=\exp\varphi$, where $\varphi:\H\to D$
is a conformal map onto a $V$-comb,
\vspace{.1in}

\noindent
(iii) $f(z)=z^qe^{-az^2+bz}\prod_{j=1}^\infty(1-z/z_j)e^{z/z_j},$
where $z_j,a,b\in\R;$
$$\sum_j|z_j|^{-2}<\infty,$$
$a\geq 0$ and $q\geq 0$ is an
integer.
\end{theorem}

It follows that there exists a function $f\in LP$ with a prescribed sequence
of critical values $c_j$ satisfying (\ref{alt}), and prescribed limits
$\lim_{x\to\pm\infty}f(x)\in\{0,\infty\}$ (asymptotic values).
Such function is defined by its critical sequence and asymptotic values
up to an increasing real affine change of the independent variable.

Here are some examples of comb representations (ii).
\begin{itemize}

\item $f(z)=z+b, b\in\R$. There are no critical values, asymptotic
values are $\pm\infty$. $D$ is a strip $\{0<|\Ima z|<\pi\}.$

\item $f(z)=\cos z$. The critical sequence is $(-1)^j$, infinite in both
directions, there are no asymptotic values. $D$ is the plane,
cut along the rays
$\{ x+i\pi j:-\infty<x\leq 0\}.$

\item $f(z)=\exp(-z^2)$. $D$ is the plane cut along the negative ray.

\item $f(z)=1/\Gamma(z)$. $D$ is the plane cut along the rays
$\{ x+i\pi j:-\infty<x\leq \log|c_j|\},\; j<0,$ where $c_j$ are
the critical values of the $\Gamma$-function, there is an asymptotic value
$0=\lim_{x\to+\infty}1/\Gamma(x)$.  

\item $f$ a polynomial of degree $n$. $D$ is obtained from the strip
$\{ x+iy:0<y<\pi n\}$ by removing the rays (\ref{rays}), where $h_j=\log|c_j|,$ and $c_j$
are the critical values of $f$.
\end{itemize}

An important subclass of $LP$ is defined
by the condition that $h_j\geq 0$, $m<j<k$,
and whenever the sequence of critical points is bounded from below (or above, or
from both sides), then the corresponding asymptotic value is $\infty$.
This subclass of $LP$ will be called the $MO$-class.
It was used for the first time in spectral theory in \cite{MO}.
Functions of $MO$-class
have another representation in terms of conformal mappings.
Consider a region $D$ of the form
\begin{equation}\label{MO-comb}
\{ x+iy: y>0,\pi m<x<\pi k\}\backslash\bigcup_{j=m+1}^{k-1}
\{\pi k+iy:0\leq y\leq h^\prime_j\},\end{equation}
where $-\infty\leq m<k\leq\infty$, and $h_j^\prime\geq 0$, see Fig. 1 (right).
Such regions will be called {\em $MO$-combs}.
Let $\theta:\H\to D$ be a conformal map $\theta(\infty)=\infty$.
Then 
\begin{equation}\label{MO}
f=\cos\theta
\end{equation}
is a function of the class $MO$ with critical values $(-1)^j\cosh h_j^\prime.$
Every function of $MO$ class can be represented in this way,
and the function is defined by its critical sequence up to a real affine
change of the independent variable.
We have the following important characterization of the $MO$ class \cite{MO}:
\begin{theorem}\label{MaO}
For a real entire function $f$,
the equation $f^2(z)-1$ has only
real roots if and only if $f\in MO$.
\end{theorem}
Such functions occur in the situation similar to the Examples 2 and 1 in
the Introduction: they describe the spectra of periodic
canonical systems \cite{Krein,deB} and entire functions of smallest
deviation from zero on closed subsets of the real axis
\cite{SY1}.

\section{MacLane's theorem}

In this section we give a geometric characterization of integrals
of $LP$ functions. Roughly speaking, we will show that critical values
of these integrals
can be arbitrarily prescribed, subject to the evident restriction
(\ref{updown}). Notice that differentiation maps $LP$ into itself,
so the class of integrals of $LP$-functions contains $LP$.

We follow the exposition in \cite{vinb} with some corrections
and simplifications, see also \cite{E1} on related questions.
Let $f$ be a real entire function
with all critical points real.

Consider the preimage $f^{-1}(\R)$.
It contains the real line, and it is a smooth curve
in a neighborhood of any
point which is not a critical point
At a critical point of order $n$ it looks
like the preimage of the real line under $z^{n+1}$.

MacLane's class consists of real entire functions for which the
preimage
of the real line looks like one of the pictures in Fig. 2,
\begin{center}
\epsfxsize=5.0in
\centerline{\epsffile{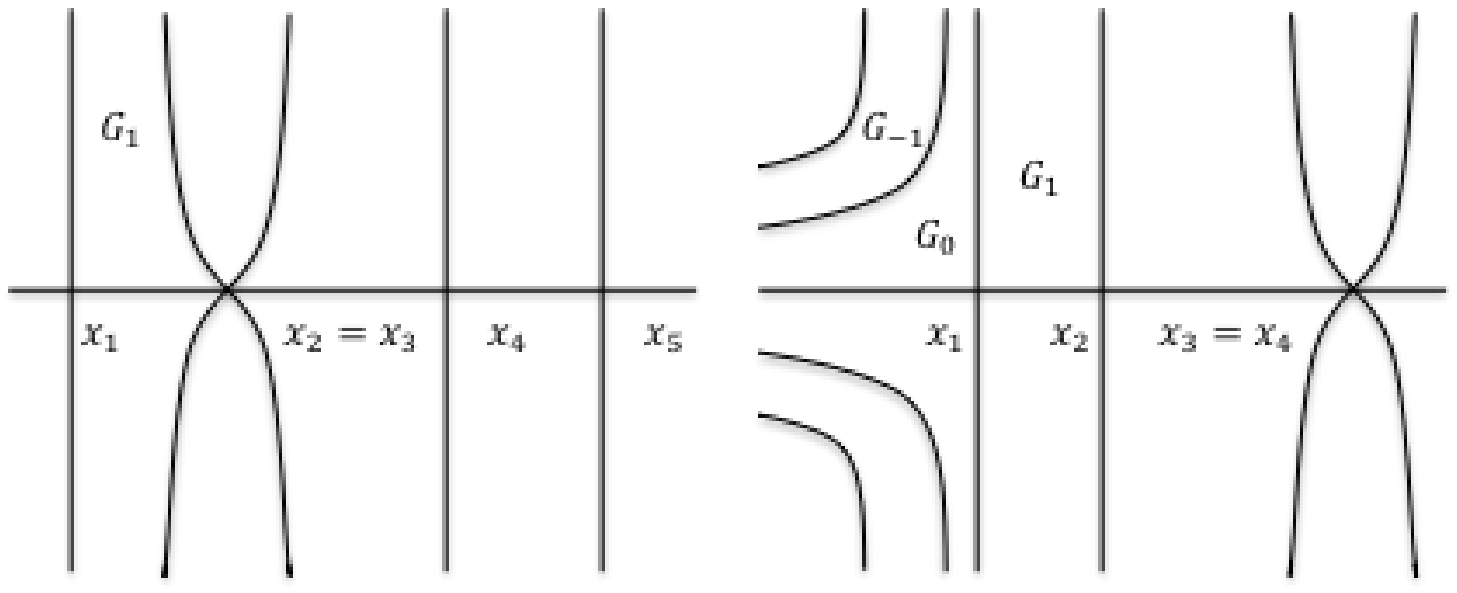}}
Fig. 2. Fish-bones.
\end{center}
up to an orientation preserving homeomorphism of the plane,
commuting with the complex conjugation. We call this picture a
fish-bone.

There are several cases. In the simplest case, the sequence
of critical points $\ldots\leq x_j\leq x_{j+1}\leq\ldots$
is unbounded from below and from above. Each critical point is
repeated in this sequence according to its multiplicity.
Preimage of the real line consists of the real line itself, crossed
by infinitely many simple curves, each curve is symmetric with
respect to the real line. The crossing points are mapped onto
the critical values $c_j=f(x_j)$. Several ``vertical''
lines cross the real line at a multiple critical point.
The complement to the union
of curves in Fig. 2, consists of simply connected regions,
each of them is mapped conformally onto the upper or lower half-plane.

The sequence of critical points can be bounded from above or
from below or both. Suppose that it is bounded from below,
and enumerate the sequence as $x_1\leq x_2\leq\ldots$.
Then the left end of the fish-bone can be
of two types. For the first type, shown in Fig 2 (left), the full preimage
of the real line is connected. We have two large complementary regions
adjacent along a negative ray, each of them is mapped
by $f$ homeomorphically onto the upper or lower half-plane.
It is easy to see that in this case we have $f(x)\to\infty$
as $x\to-\infty$. We set $c_0=\infty$ and extend our critical
sequence $(c_j)$ by adding this term to it.

The second type of the end is shown in Fig. 2 (right).
In this case, the preimage
of the real has infinitely many components. In addition to one component
of $f^{-1}(\R)$, as above, there are infinitely many simple curves tending to infinity
at both ends. Strip-like regions between these curves are mapped
homeomorphically onto the upper or lower half-plane.
In this case $c_0=\lim_{x\to-\infty} f(x)\neq\infty$, and we extend our critical sequence by
$c_0$.

Similar situations may occur on the right end when the sequence
of critical points is bounded from above.

In all cases, the fish-bone is completely determined by the
{\em augmented critical
sequence} $(c_j)$. We use the following notation: if the
sequence of critical points is unbounded from above and below, then
$-\infty<j<+\infty$. In all other cases, the critical values are
$c_j;\; m<j<k$, where $-\infty\leq m<k\leq+\infty$, and if $m$ or $k$
or both are finite, we add to our sequence the term $c_m$ or $c_k$
or both, which are the limits of $f(x)$ as $x\to-\infty$ or
$x\to+\infty$.

The augmented critical sequence 
satisfies the following condition
\begin{equation}\label{updown}
(c_{j+1}-c_j)(c_{j}-c_{j-1})\leq 0.
\end{equation}
All $c_j$ are real, except possibly the first and/or the last
term which can be $\pm\infty$.

We call such sequences ``up-down sequences''. If the sequence of
critical points is unbounded from below and from above, then 
the sequences $x_j$ and $c_j$ are defined for a given $f$
up to a shift of the
subscript.
\vspace{.1in}

\noindent
{\bf MacLane's theorem} \cite{Mc} {\em For every up-down sequence,
finite or infinite
in one or both directions, there exists a function
$f\in M$ for which this
sequence is the critical sequence.
Any two functions corresponding to the
same sequence are related by $f_1(z)=f_2(az+b)$ with $a>0,b\in\R$.} 
\vspace{.1in}

In other words, one can prescribe a piecewise-monotone
graph on the real line, and after
a strictly increasing continuous change of the independent variable,
this will be
the graph of an entire function of MacLane's class, which is essentially
unique.

Uniqueness statement in MacLane's theorem is easy.
Suppose that $f_1$ and $f_2$
are two functions of MacLane class with the same augmented critical
sequence. Then it is easy to construct a homeomorphism 
$\phi$ of the plane such that $f_1=f_1\circ\phi$. Then $\phi$
must be conformal and commute with complex conjugation,
so $\phi$ must be a real affine map.

Class $LP$ is contained in the MacLane class.
It corresponds to the case when the critical sequence
satisfies the condition (\ref{alt})
and in addition, the first and last terms of the sequence $c_j$, if present, are
$0$ or $\infty$. It is clear that (\ref{alt})
is stronger than (\ref{updown}).

We proved this special case of MacLane's theorem in the previous
section.
Now we give the proof of MacLane's theorem in full generality.
\vspace{.1in}

First we recover the fish-bone from the given sequence $(c_j)$ as explained above.
Then we construct a continuous map $F:\C_z\to\C_w$ as follows. We map each interval
$[x_j,x_{j+1}]\in\R$ linearly onto the interval $[c_j,c_{j+1}]$.
Then we map each infinite ray of the fish-bone onto a corresponding ray 
of the real line, linearly with respect to length. The curves on the left
of Fig. 2 (right) are mapped on the rays $[c_0,\infty)$.
Then we extend our map to the components of the complement of the fish-bone, so that
each component is mapped on the upper or lower half-plane homeomorphically.

The resulting continuous map $F$ is a local homeomorphism everywhere except the
points $x_j$ where it is ramified. There is unique conformal structure $\rho$
in the plane $\C_z$ which makes this map holomorphic. By the uniformization theorem,
the simply connected Riemann surface $(\C,\rho)$ is conformally equivalent
to a disc $|z|<R$, where $R\leq\infty$. This means that there exists a homeomorphism
$\phi:\{ z:|z|<R\}\to\C$ such that $F\circ\phi$ is a holomorphic function.
As all our construction
can be performed symmetrically with respect to the real line, $F$ is a real function.
It remains to prove that $R=\infty$.

If the sequence $(x_j)$ is finite, and both asymptotic values are $\infty$,
our map extends to a continuous map
of the Riemann sphere $\bC_z\to\bC_w$ by putting $F(\infty)=\infty$.
So the Riemann surface $(\bC,\rho)$ must be conformally equivalent to the sphere,
and we obtain that $R=\infty$. In this case $f=F\circ\phi$
is evidently a polynomial.

If the sequence $(x_j)$ is infinite in both directions, we consider truncated sequences
$(c_j)_{j=-n}^n$, augmented by asymptotic values $\infty$ on both sides,
and the corresponding fish-bones and maps $F_n$ as above. By the previous argument 
we have homeomorphisms $\phi_n$ and polynomials $f_n=F_n\circ\phi_n$.

We can always arrange that $x_1<x_2$, $0\in(x_1,x_2)$, and $F_n(0)=a\in(c_1,c_2),$
where $a$ is independent of $n$. Then we choose $\phi_n$ so that $\phi_n(0)=0$,
and 
\begin{equation}\label{norm}
f_n^\prime(0)=1.\end{equation}
Then $f_n$ maps univalently some disc $\{ z:|z|<r\}$ onto a region $G_n$ which contains
a disc $\{ w:|w|<\epsilon\}$ and is contained in a disc $\{ w:|w|<\delta\}$
with some $r>0,\epsilon>0,\delta>0$ which are independent of $n$. This follows from the
Schwarz lemma applied to $f_n$ and $f_n^{-1}$ in a neighborhood of $0$.
We conclude that $(f_n)$ is a normal family in $\{ z:|z|<r\}$ and the limit functions
are non-constant.

Now we use the following lemma \cite{L1}.

\begin{lemma}\label{lemma2} Let $g_n$ be a sequence of
real polynomials whose all zeros are real,
and suppose that $g_n\to g\not\equiv 0$ uniformly in some
neighborhood of $0$. Then
$g$ is entire, and $g_n\to g$ uniformly on compact subsets of $\C$.
\end{lemma}

{\em Proof.} By a shift of the independent variable we may assume that
$g(0)\neq 0.$ Then $g_n(0)\neq 0$ for large $n$. We have 
$$-\left(\frac{g_n^\prime}{g_n}\right)'(0)=\sum_k\frac{1}{z_{n,k}^2},$$
where $z_{n,k}$ are zeros of $g_n$.
The left hand side is bounded by a constant independent of $n$, while
all summands in the right hand side are positive. So for every interval $I$
on the real line there exists a constant $c(I)$ independent of $n$ such that
the $g_n$ have at most $c(I)$ roots on $I$. Thus from every sequence of $g_n$
one can choose a subsequence such that the zero-sets of polynomials
of this subsequence tend to a limit set which has
no accumulation points in $\C$. So our subsequence converges to an
entire function. Evidently this entire function is an analytic
continuation of $g$, and the statement of the lemma follows.
\vspace{.1in}

We apply this lemma to the sequence $(f_n^\prime)$ and conclude that $f$ is entire, that is
$R=\infty$, as advertised.
\vspace{.1in}

Now we describe the necessary modifications of this proof for the case that the sequence of critical points is bounded
from below (the case of semi-infinite sequence bounded from above is treated similarly).
If the asymptotic value $c_0=\infty$, no modification is needed.
If $c_0\neq \infty$, we may assume without loss of generality that $c_0=0$, by adding a real constant
to all functions $f,F,f_n$. Then we approximate our critical sequence $c_0,c_1\ldots$ by the finite sequences
$c_0,c_0,\ldots,c_0,c_1\ldots,c_n$, where $c_0=0$ is repeated $n$ times. The corresponding fish-bone is shown
in Fig. 3, where $\beta$ is the additional zero of multiplicity $n$.
As $n\to\infty$, $\beta\to-\infty.$
 The rest of the argument goes without change.
\begin{center}
\epsfxsize=4.0in
\centerline{\epsffile{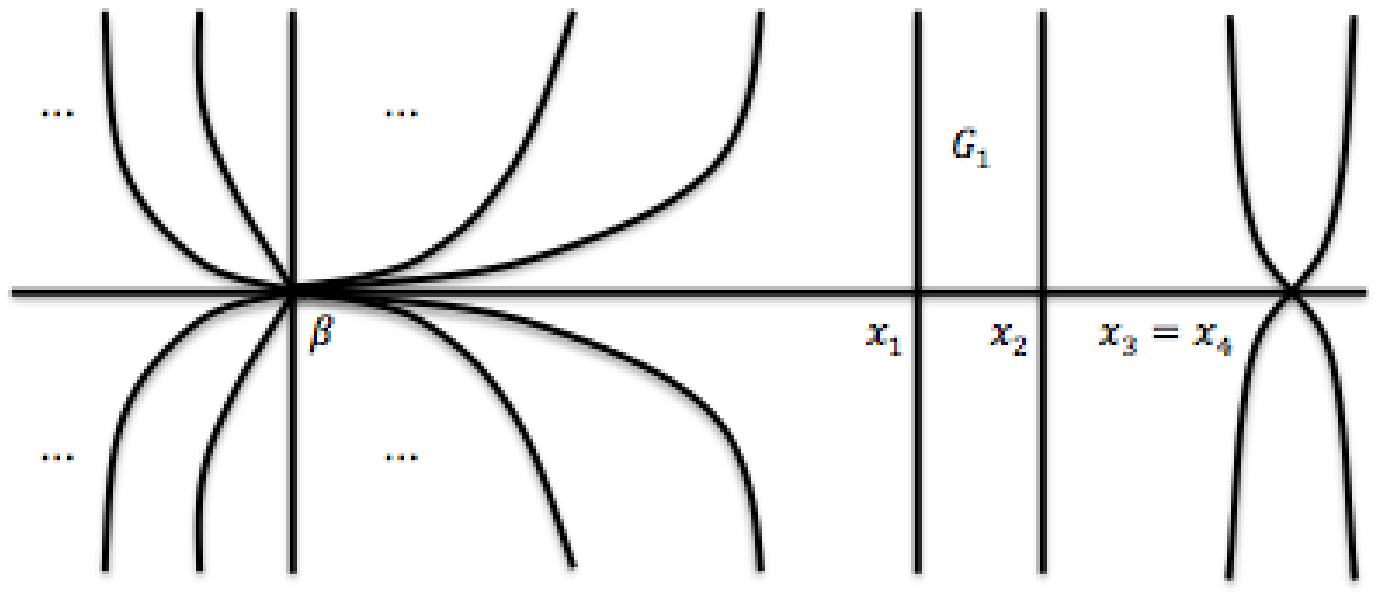}}
Fig. 3. Approximation of a fish-bone by polynomial ones.
\end{center}
Finally we consider the case when there are finitely many critical points and two different asymptotic values.
In this case, no approximation argument is needed, and $f(z)=\int_{-\infty}^z P(\zeta)\exp(-a\zeta^2+b\zeta)d\zeta$, where
$P$ is a real polynomial with all zeros real and $a\geq 0$ and $b\in\R$.

\section{Representation of Green's and Martin's\break functions}

Here we discuss the relation between comb regions
and Green and Martin functions of complements
of closed sets on the real line.

Let $E\subset\R$ be a compact set of positive capacity.
Then there exists the Green function $G$
of $\Omega=\bC\backslash E$ with pole at $\infty$.
We have
\begin{equation}\label{pot}
G(z)=\int_E\log|z-t|d\mu(t)+\gamma(E),
\end{equation}
where $\mu$ is a probability measure on $E$ which is called the equilibrium
measure, and $\gamma$ the Robin constant of $E$. 
Function $G$ is positive and harmonic in $\C\backslash E$, and has boundary
values $0$ a. e. with respect to $\mu$. We have
\begin{equation}\label{robin}
G(z)=\log|z|+\gamma+o(1),\quad z\to\infty.
\end{equation}
These properties characterize
$G$ and $\mu$ \cite{Land}.

There exists an analytic function $\phi:\H\to \H$, such that $G=\Ima\phi$.
It is called the complex Green function. Since the derivative
$$\phi'=\frac{d}{dz}\left(i\int_E\log(z-t)d\mu(t)\right)=
i\int_E\frac{d\mu(t)}{z-\zeta}$$ 
has positive real part in $\H$, we conclude from Lemma~\ref{lemma1}
that $\phi$ is univalent. Let $D=\phi(\H)$. This region $D$ has the following
characteristic properties:
\vspace{.1in}

\noindent
(i) $D$ is contained in a vertical half-strip $\{ x+iy:a<x<b,\; y>0\}$
with $b-a=\pi$,
and contains a half-strip $\{ x+iy:a<x<b,\; y>K\}$ with some $K>0$.
\vspace{.1in}

\noindent
(ii) For every $z\in D$, the vertical ray $\ell_z=\{ z+it:t\geq 0\}$ is
contained in $D$.
\vspace{.1in}

\noindent
(iii) For almost every $x\in(a,b)$, the ray $\{ x+iy:y>0\}$ is contained in
$D$.
\vspace{.1in}

These properties can be restated shortly as follows:
\begin{equation}\label{h}
D=\{ x+iy:a<x<b,\; y>h(x)\},
\end{equation}
where $h$ is as non-negative upper semi-continuous function bounded
from above and equal to $0$ a. e. 
\vspace{.1in}

\noindent
{\em Remarks.} Function $G$ given by (\ref{pot}) is upper semi-continuous,
so it must be continuous at every point where $G(z)=0$.
If $h(x)=0$ for some $x\in(a,b)$, then for the similar reason,
$h$ is continuous at $x$, so $\partial D$ is locally connected at $x$.
It follows that $x=\phi(x')$ for some $x'\in\R$, and $\phi$ is continuous at
$x'$. In other words, existence of a radial limit $\phi(x')$, such that
$\Ima\phi(x')=0$ implies continuity of $\phi$ and $G$ at the point $x'$.
\vspace{.1in}

To prove (i) we integrate by parts: 
\begin{equation}\label{poiss}
\Rea\phi(z)=-\int_E\arg(z-t)d\mu(t)=y\int_E\frac{\mu(t)dt}{(x-t)^2+y^2},
\end{equation}
where $z=x+iy$ and $\mu(t)=\mu((-\infty,t])$ is the distribution
function. As $\mu$ has no atoms, $t\mapsto\mu(t)$ is continuous.
So $\Rea\phi$ is continuous in $\overline{\H}$.
The first statement of (i) follows because $0\leq \mu(t)\leq 1$, and
the second because $G(z)=\Ima\phi$ is bounded on any compact set in $\C$.

To prove (ii), let $\alpha$ be a tangent vector to the ray $\ell_z$,
so $\alpha=i$.
Then $\beta=(\phi^{-1})' \alpha$
will be the tangent vector to the $\phi$-preimage
of this ray, and we have seen that $\arg(\phi^{-1})'\in(-\pi/2,\pi/2)$.
So $\beta$ is in the upper half-plane thus the preimage of $\ell_z$ 
can never hit the real line, and an analytic continuation of $\phi^{-1}$
is possible along the whole ray $\ell_z$.

To prove (iii), we use (\ref{poiss}) again. 
As $\mu(t)$ is continuous, $\Rea\phi$ is continuous in $\overline{\H}$.
Moreover, 
$$\Rea\phi(\beta)-\Rea\phi(\alpha)=\mu(\beta)-\mu(\alpha),\quad \alpha<\beta.$$
This means that measure $\mu$ on $E$ corresponds to the Lebesgue measure on
base of the comb $(a,b)$. Furthermore, if for some $x\in E$ we have
$G(x)=0$ then $h(\Rea\phi(x))=0$. Thus $h=0$ almost everywhere with respect to
the Lebesgue measure on $(a,b)$. This proves (iii),

Now we show that for every $D$ satisfying (i)--(iii), the
conformal map $\phi:\H\to D$ is related with the Green function $G$
of some closed set $E$ by the formula $G=\Ima\phi$.

Imaginary part $v=
\Ima \phi$ is a positive harmonic function  in the upper half-plane.
We extend it
to the lower half-plane by symmetry, $v(\overline{z})=v(z)$,
and to the real line by upper semicontinuity: $v(x)=\limsup_{z\to x}v(z)$.
In view of (i), $\partial D$ has a rectilinear part near infinity,
the extended function $v$ is harmonic in
a punctured neighborhood of $\infty$ and has asymptotics of the form
$$v(z)=\log|z|+\mathrm{const}+o(1),\quad z\to\infty.$$
Let us prove that $v$ is subharmonic in the whole plane, and has a
representation
(\ref{pot}) with some probability measure $\mu$ 
with compact support on the real line.

Let $\{ h_k\}$ be a dense set on $\partial D$.
Let $D_n$ be the region obtained from the half-strip
$\{ x+iy:a<x<b,\; y>0\}$ by removing the vertical segments $\{\Rea h_k+iy,\; 0<y\le\Ima h_k\}$.
Then $D_1\supset D_2\supset\ldots\to D$. Let $\phi_n$ be conformal maps
of $\H$ onto $D_n$, normalized by $\phi_n(0)=a,\;\phi_n(1)=b,\;
\phi_n(\infty)=\infty$. Then it is easy to check that
$\Ima\phi_n$ is the Green function of
some set $E_n\subset[0,1]$ consisting of
finitely many closed intervals.
So 
$$\Ima\phi_n(z)=\int\log|z-t|d\mu_n(t)+\gamma_n,$$
with some probability measures $\mu_n$ on $[0,1]$ and some constants $\gamma_n$.
We can choose a subsequence 
such that $\mu_n\to\mu$ weakly, where $\mu$ is a probability measure on
$[0,1]$, and it is easy to check that (\ref{pot}) holds with some $\gamma$.
Thus $v$ is subharmonic
in the plane. Since $v\geq 0$, the measure $\mu$ has
no atoms. It remains to prove that $v(x)=0$ a. e. with respect to
$\mu$. This follows from the property (iii) of the region $D$.
Indeed, let $x\in(a,b)$ be a point such
that the vertical ray $\ell_x$ is in $D$,
except the endpoint $x$. By a well-known argument,
the curve $\phi(\ell_x)$ has an endpoint at some $x'\in(0,1)$,
and the angular limit
of $v=\Ima\phi$ is zero at this point $x'$. By the remark above, $v(x')=0$.
We define $E$ as the closed support of $\mu$. Then $v=\Ima\phi$ is
positive and
harmonic outside $E$ and $v(x)=0$ $\mu$-almost everywhere,
so $v$ is the Green
function of~$E$.

Our construction of $D$ from $E$ defines $D$ up to a shift by a real number.
The inverse construction defines $E$ up to a real affine transformation,
and changing $E$ by a set of zero capacity.

Now we give a similar representation of Martin functions.
Let $E\subset\R$ be an unbounded closed set of positive capacity. Let
$U$ be the cone of positive harmonic functions in $\C\backslash E$,
and $U_s\subset U$ the cone of symmetric positive harmonic functions,
$v(z)=v(\overline{z})$.
Martin's functions are minimal elements of $U$, that is functions
$v\in U$ with the property
$$u\in U,\quad u\leq v\quad\mbox{implies}\quad u=cv,$$
where $c>0$ is a constant. 
Similarly we define {\em symmetric Martin functions} using $U_s$ instead
of $U$. Martin functions always exist and form a convex cone.
If $v$ is a Martin function, then
$u(z)+u(\overline{z})$ is a symmetric Martin function,
so symmetric Martin functions
also exist and form a convex cone.

Let $v$ be a symmetric Martin function, and let $\phi$ be an analytic function
in $\H$ so that $v(z)=\Ima\phi(z),\; z\in \H$. Then $\phi:\H\to D$ is a conformal
map onto a region $D\subset \H$. This is proved in the same way as for
Green's functions. Regions $D$ arising from symmetric Martin functions are
characterized by the properties (ii), (iii) above and the {\em negation} of
the property (i): either $a=-\infty$ or $b=+\infty$, or $h$ is unbounded in
(\ref{h}). 

Notice that function $\phi$ maps $h$ into $\H$, so the angular derivative
of $\phi$ at infinity exists, that is
$$\phi(z)=cz+o(z),\quad z\to\infty\quad\mbox{in any Stolz angle},$$
where $c\geq 0$. One can derive from this that every Martin function
satisfies 
$$B(r,v):=\max_{|z|=r}v(z)=O(r),\quad r\to\infty.$$
This implies that the cone of Martin functions has dimension at most $2$,
\cite{Kjellberg,L5,deB}, and the cone of symmetric Martin functions is always
one-dimensional.

Dimension of the cone of Martin's functions is an important characteristic
of the set $E$, see \cite{Benedicks,L5}. One can show that
the cone of Martin functions is two-dimensional if and only if
$$\limsup_{r\to\infty}B(r,v)/r>0.$$

Now we impose various conditions on $E$ and find their exact counterparts
in terms of $E$ and $\mu$.

The first important condition is that the set $E$ is regular
in the sense of potential theory
\cite{Land}. In this case Green's and Martin's functions are continuous in
$\C$. For the region $D$ this is equivalent to the local connectedness
of $\partial D$ in the case of Green's function, and local connectedness
of the part $\partial D\backslash X$,
where $X$ is the union of the vertical rays
on $\partial D$, if these rays are present.
In terms of function $h$ in (\ref{h}),
local connectedness is equivalent in the case of Green's function 
to the
condition that the set $X=\{ x: h(x)>0\}$ is at most countable, and the sets
$X_\epsilon=\{ x:h(x)>\epsilon\}$ are finite for every $\epsilon>0$,
that is $D$ is obtained from $\H$ by making countably many cuts, and the length
of a cut tends to $0$.

In the case of Martin's function, local connectedness of $D$ means that
the sets $X_\epsilon$ can only accumulate to $a$ or $b$.

Next we discuss the condition on $D$ which corresponds
to absolute continuity of $\mu$. We thank Misha Sodin who passed to us the
contents of his conversation with Ch. Pommerenke on this subject.

To state the result we first recall McMillan's sector theorem \cite{McM},
\cite[Thm. 6.24]{Pom}.
Let $f$ be a conformal map from $\H$ to a region $G$.
Let $\sect(f)$ be the set of points $x\in\R$ such that the non tangential
limit $f(x)$ exists and $f(x)$ is the vertex of an angular sector in $G$.
\vspace{.1in}

\noindent
{\bf McMillan Sector Theorem} \cite[Theorem 6.24, p.146]{Pom}.
{\em 
Assume that $A\subset \sect(f)$. Then
\begin{equation}
\label{mcmillan}
|A|=0 \quad\text{if and only if}\quad |f(A)|=0.
\end{equation}
}

We say that the sector condition holds in the comb region $D$ if the function
\begin{equation}
\label{sector}
H(x)=\sup_{y\in (a,b)}\frac{h(y)}{|y-x|}
\end{equation}
is finite for almost all $x\in (a,b)$.
Geometrically it means that for almost all
$x$ in the base of the comb there exists a Stolz
angle with the vertex in $x$. 

\begin{theorem} \label{th4} Region $D$ satisfies the sector
condition if and only if
$\mu$ is absolutely continuous with respect to the Lebesgue measure on $\R$.
\end{theorem} 

{\em Proof}.
Recall that the Lebesgue measure on the base of the comb
corresponds to the harmonic measure $\mu$ on $E$.

Assume that the sector condition holds.
This means that a Borel support of the harmonic measure
$\mu$ is contained in $\sect(\phi)$.
Let $A$ be a Borel support of the singular component of $\mu$.
By the definition $|A|=0$. Thus, by 
McMillan's theorem $\mu(A)=0$,
thus $\mu$ is absolutely continuous.

Conversely, assume that the harmonic measure is absolutely continuous.
Recall that $\phi'$ has positive imaginary part,
and therefore possesses non-tangential
limits for almost all $x$ with respect to the Lebesgue measure.
Therefore the limit exists for almost all
$x$ with respect to the harmonic measure as well.
\medskip

\noindent
{\em Example 1.}  There exist irregular regions with absolutely continuous
measures $\mu$. Indeed, let $C$ be the standard Cantor set in $[a,b]$.
Let $h(x)$ be the characteristic function of $C$.
Then the region generated by this comb is irregular,
on other hand $H(x)$ is finite for all $x\in[a,b]\setminus C$.
\medskip

\noindent
{\em Example 2.} We give an example of a
comb such that the conditions of the previous theorem do not hold,
moreover $H(x)=\infty$ for almost all $x\in [a,b]$.
This comb is related to the Julia set
of a polynomial $T(z)=z^2-\lambda$ \cite{SY3}.
For $\lambda>2$ there exists $h_0>0$ such that the Julia
set of $T$ is the preimage of the base of the comb given in Fig. 4.
\begin{center}
\epsfxsize=3.5in
\centerline{\epsffile{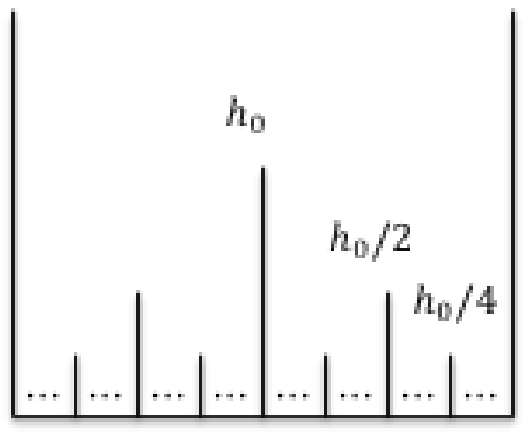}}
Fig. 4.
Comb related to the Julia set of $T(z)=z^2-\lambda$.
\end{center}
Recall that almost every number $x$ contains arbitrarily long 
strings of zeros in its dyadic representation, that is, for almost every $x$,
and every non-negative integer $N$, there exists a string $y_m$ of $0$'s and
$1$'s, ending with $1$, such that
$x=(y_m,\overbrace{0,...,0}^N,....)$.
Then $h(y_m)=2^{-m}h_0$,
and $|x-y_m|=2^{-(m+N)}$.
That is $H(x)\ge 2^Nh_0$.
In fact the Lebesgue measure of the Julia set is $0$,
i.e. the harmonic measure is singular continuous.
Note that since the boundary is locally connected the region
$\bC\setminus E$ is regular.
\medskip
 
Even stronger condition is that 
\begin{equation}\label{widom}
\sum_x h(x)<\infty,
\end{equation}
in other words, the total length of slits is finite.
This is the so-called Widom condition.
It appears in his studies of asymptotics for 
extremal polynomials associated with a system of curves in the complex plane.
Let $\pi_1(\Omega)$ be the fundamental group of the given region
$\Omega=\overline{\C}\setminus E$. 
To a fixed character $\alpha\in \pi^*_1(\Omega)$
one associates the set of multi-valued 
(character-automorphic) analytic functions
$$
H^\infty(\alpha)=\{f:f\circ\gamma=
\alpha(\gamma)f,\;\forall\gamma\in \pi_1(\Omega),
\; \sup_{z\in\Omega}|f(z)|<\infty\}.
$$
The region $\Omega$ is of Widom type if
the space $H^\infty(\alpha)$ is non-trivial 
(contains a non-constant function) for every $\alpha\in \pi^*_1(\Omega)$. 
A regular region $\Omega=\overline{\C}\setminus E$
is of Widom type if and only if (\ref{widom}) holds.
For the role of this condition in the spectral
theory of almost periodic Jacobi matrices see \cite{SY2}.
A well-known fact that the derivative of a conformal mapping
on a region bounded by a rectifiable curve belongs to $H^1$
implies that the corresponding equilibrium measure $\mu$
is absolutely continuous.

\section{More general combs}

In this section we consider more general comb regions:
those which satisfy property (ii) of the previous section.
These regions $D$ can be described as
\begin{equation}\label{usc}
D=\{ x+iy:a<x<b,\; y>h(x)\},
\end{equation}
where $-\infty\leq a<b\leq\infty$ and $h$ is an upper semi-continuous function
on $(a,b)$.
Let $\Combs$ be the set of such regions modulo horizontal shift, or
equivalently, the set of all triples $(a,b,h)$ modulo the equivalence relation
$(a,b,h)\sim (a+c,b+c,h(x-c))$, $c\in\R$.

Let $\Conf$ be the set of univalent functions in $\H$ such that
$\phi(\H)\in\Combs$, normalized by $\phi(\infty)=\infty$, modulo the
equivalence relation $\phi(z)\sim\phi(z-c),\; c\in\R$. 

Let $R$ be the set of all analytic functions with positive real
part in $\H$.

Let $\Sub$ be the set of all subharmonic functions $v$ in the plane
of the form
$$v(z)=\Rea\left(-az^2+bz+\int_\R\left(\log\left(1-\frac{z}{t}\right)+
\frac{zt}{1+t^2}\right)d\mu(t)\right),$$
where $a>0,\; b\in\R$ and
$\mu$ is an increasing right-continuous function,
such that
$$\int_{0}^\infty\frac{\mu(t)-\mu(-t)}{1+t^2}dt<\infty.$$
Two such functions are considered equivalent if their difference is
constant.

\begin{theorem} There are the following canonical bijections between
the sets $\Conf$, $R$, $\Sub$:

$$\phi\mapsto \phi':\Conf\to R,\quad
\phi\mapsto \Ima \phi:\quad\Conf\to\Sub.$$

Moreover, $\Rea \phi=\mu+\const$, $(2\pi)^{-1}\Delta v=d\mu$.
\end{theorem}

\section{Uniform approximation}

Here we consider several extremal problems whose solutions are
expressed in
terms of comb functions.

Applications of comb functions to extremal problems begins with the
work of Akhiezer and Levin \cite{AL} on extension of Bernstein's inequality.
Further applications are contained in \cite{L2,L3,L4,L5}. A survey of
polynomials and entire functions of least deviation from zero on
closed sets on the real line is
given in \cite{SY1}.

Here we mention only few results.
\vspace{.1in}

1. Let $f$ be an entire function of exponential type $1$ satisfying
$|f(x)|\leq 1,\; x<0$ and $|f(x)|\leq B,\; x>0,$ where $B\geq 1$. 
One looks for maximal values of $|f(x)|$ for given $x$ and of $|f'(0)|$,
\cite{E}.
The extremal function is expressed in terms of the $MO$-comb
with $h_j^\prime=0,\; j< 0$ and $h_j^\prime=
\cosh^{-1} B,\; j\geq 0.$ Let $\theta:\H\to D$ be the conformal
map onto the region (\ref{MO-comb}), such that $\theta(z)\sim z,$ as
$z\to\infty$ non-tangentially,
$\theta(0)=\pi-$. Set $x_1=\theta^{-1}(ih_0)$.
Then the function
$$f_0(x)=\left\{\begin{array}{ll}B,&x>x_1,\\
                              \cos\theta(x),& 0\leq x\leq x_1,\\
                              1,&x<0
\end{array}\right.$$
gives the solution of the first extremal problem: $|f(x)|\leq f_0(x)$
for $f$ in the class described above, and $f_0^\prime(0)$ is the maximal
value of $|f'(0)|$.
\vspace{.1in}

2. Best uniform approximation of $\sgn(x)$ on two rays/intervals.
The simplest problem of this kind is to find the best uniform approximation
of $\sgn(x)$ on the set $X=(-\infty,-a]\cup[a,\infty)$
by entire functions of exponential type at most $1$. The 
extremal entire function
belongs to the MacLane class and has critical sequence
\begin{equation}
c_j=\left\{\begin{array}{ll}-1+(-1)^jL,& j\leq 0,\\
                                1+(-1)^jL,&j>0,\end{array}\right.
\end{equation}
where $L$ is the error of the best approximation.
Unfortunately, MacLane's functions do not have simple representations
in terms of conformal mappings
like (\ref{vc})
or (\ref{MO}), however in certain cases representation in terms of
conformal maps of the kind described in section 5
can be obtained \cite{EY1,EY3}.
\vspace{.1in}

\def\cJ{\mathcal{J}}

3. Let us consider a uniform counterpart of the classical orthogonal Jacobi polynomials.
Let $\alpha,\beta\ge 0$ and let $\cJ_n(x;\alpha,\beta)=x^n+\dots$
denote the monic polynomial of least deviation from zero on $ [0,1]$
with respect to the weight function $x^\alpha(1-x)^\beta$.
\begin{lemma}\label{lemma3} For non-negative  $\alpha,\beta$ and an integer $n$
$$x^\alpha(1-x)^\beta\cJ_n(x)=L e^\phi,$$
where $\phi$ is the conformal map on the $V$-comb region
$$
D=\{z=x+iy:-\beta<\frac y \pi<\alpha+n\}\setminus
\bigcup_{j=0}^n\{z=x+iy:\frac y\pi=j,\ x\le 0\}.
$$
\end{lemma}
Such polynomials turn out to be useful in the
description of multidimensional 
polynomials of least deviation from zero \cite{moyu}.
As an example we formulate the following theorem.
Note that in multidimensional situation an extremal polynomial
is not necessarily unique.

\begin{theorem} $\cite{mo}$ A best polynomial approximation $P(z_1,\dots,z_d,
\overline{z_1},\dots. \overline{z_d})$   to the monomial $z_1^{k_1}\dots z_d^{k_d}\overline{z_1}^{l_1}$, $k_1\ge l_1$, by polynomials of the total degree less than
$k_1+\dots+k_d+l_1$ in the ball $|z_1|^2+\dots+|z_d|^2\le 1$ can be given in the form
\begin{equation*}
\begin{split}
z_1^{k_1}\dots z^{k_d}\overline{z_1}^{l_1}+P(z_1,\dots,z_d,
\overline{z_1},\dots. \overline{z_d})\\
=
z_1^{k_1-\ell_1}z_2^{k_2}\dots z_d^{k_d}\cJ_{l_1}\left(|z_1|^2;\frac{k_1-l_1}{2}, \frac{k_2+\dots+k_d} 2\right).
\end{split}
\end{equation*}

\end{theorem}
\section{Spectral theory and harmonic analysis}

1. We say that an unbounded closed set $E$ is homogeneous if there exists $\eta>0$ such that for all $x\in E$ and all $\delta>0$,
$$|(x-\delta,x+\delta)\cap E|\ge \eta\delta.$$
\begin{theorem}$\cite{Y}$
Let $\theta$ be a conformal map from the upper half-plane $\H$ onto
an $MO$-comb region $D$ (Fig. 1, right).
Assume that $E=\theta^{-1}(\R)$ is homogeneous.
Then $E$ is the spectrum of a periodic canonical system, i.e., there exists
an integrable on $[0,1]$ non negative $2\times 2$ matrix function $H(t)$  of period 1, $H(t+1)=H(t)$, such that for
an entire (transfer) matrix function $T(1,z)$ defined by the differential system
\begin{equation}\label{can1}
J\dot T(t,z)=zH(t)T(t,z),\ \ T(0,z)=I,
\quad J=\begin{bmatrix}0&1\\-1& 0\end{bmatrix},
\end{equation}
 the following relation holds 
\begin{equation}e^{i\theta}=
\Delta-\sqrt{\Delta^2-1},\quad \Delta(z):=(1/2)\tr T(1,z).
\end{equation}
Moreover the parameter $t$ in \eqref{can1}
corresponds to the ``exponential type'' of the matrix $T(t,z)$
with respect to the Martin function $\theta$, that is,
\begin{equation}
t=\lim_{y\to +\infty}\frac{\log\|T(t,iy)\|}{\Ima \theta(iy)}.
\end{equation}
The whole collection of such matrices $H(t)$
for the given $E$ can be pa\-ra\-metri\-zed by
the characters of the funda\-men\-tal group of the region
$\Omega=\C\setminus E$.
\end{theorem}

The condition of homogeneity of $E$ implies that $\Omega=\C\backslash E$
is of Widom type,
and thus $D$ satisfies Widom's condition (\ref{widom}).
This fact plays a crucial tole in the proof of Theorem 7.

\noindent {\em Example.} A region $D$ is defined by a system of slits forming a geometric progression\begin{equation*}
h_{j^k}=\kappa j^k, \ \kappa>0,\ \ h_0=\infty,
\end{equation*}
otherwise $h_j=0$. The corresponding set $E$ is homogeneous.
\vspace{.1in}

2. Riesz bases. A sequence of vectors $(e_n)$ in a Hilbert space $H$ is
called a {\em Riesz basis} if it is complete and there exist positive constants $c,C$ such that
$$c\sum|a_n|^2\leq\left\|\sum a_ne_n\right\|\leq C\sum|a_n|^2$$
for every finite sequence $(a_n)$. A long-standing problem is how to find out
whether for a given sequence of real exponents $(\lambda_n)$ the sequence
$e^{i\lambda_nx}$ is a Riesz basis in $L^2(-\pi,\pi)$.
A recent result of Semmler gives a parametric description of such Riesz bases.

We say that a sequence $(d_n),\; d_n\geq 0$ satisfies the
{\em discrete Muckenhoupt
condition} if 
$$\sum_{n\in I}d_n\sum_{n\in I}d_n^{-1}\leq C(\card I)^2,$$
for every interval $I$ of integers, and some $C>0$.

\begin{theorem} $\mathrm{\cite{Sem}}$ The sequence $(e^{i\lambda_n x})$ is a Riesz 
basis in $L^2(-\pi,\pi)$ if and only if it is the sequence of zeros of
the entire function $f=\exp\phi$ of exponential type,
where $\phi$ is a conformal map
onto a $V$-comb with tips of the cuts $h_n$, and $\exp(2h_n)$ satisfies
the discrete Muckenhoupt condition.

For a given sequence $(h_n)$ such that $(\exp(2h_n))$ satisfies the discrete
Muckenhoupt condition, the conformal map $\phi$ can be always normalized
so that $f=\exp\phi$ is of exponential type.
\end{theorem}

This theorem parametrizes all Riesz bases consisting
of functions $e^{i\lambda_nx}$ in terms of sequences $h_k$.
\vspace{.1in}

We thank Misha Sodin for many illuminating discussions on the subject
in the period 1980--2011.

\vspace{.1in}

{\em Department of Mathematics, Purdue University

West Lafayette, IN 47907 USA

eremenko@math.purdue.edu
\vspace{.1in}

Abteilung f\"ur Dynamische Systeme und Approximationstheorie,

Johannes Kepler Universit\"at Linz,

A--4040 Linz, Austria

Petro.Yudytskiy@jku.at}
\end{document}